\documentclass[preprint,fleqn,11pt]{elsarticle}
\usepackage{amssymb}
\biboptions{comma,square}

\journal{Journal of Computational and Applied Mathematics}
\usepackage{amssymb,amsmath}
\usepackage{graphics}
\usepackage{graphicx}
\usepackage{color}

\newcommand{\dl}[1]{{\bf Theorem{#1.}}}
\newcommand{\yl}[1]{{\bf Lemma{#1.}}}

\newcommand{\zb}{\textbf{\qquad$\Box$}}
\newcommand{\la}[1]{\label{#1}}
\newcommand{\rf}[1]{(\ref{#1})}

\newcommand{\bo}{\boldsymbol}

\newcommand{\zm}{{\bf Proof.}}
\newcommand{\commentout}[1]{{}}

\begin{document}
\thispagestyle{empty} \setcounter{page}{1}

\begin{frontmatter}

\title{The Weak Galerkin Finite Element Method for the Symmetric Hyperbolic Systems}

\author{Tie Zhang\footnote{Corresponding author at: Department of
Mathematics, Northeastern University, Shenyang 110004, China. {\em
E-mail address} : ztmath@163.com (T. Zhang). Tel \& Fax:
+86-024-83680949. },\quad Shangyou Zhang$^\dag$
}
\address{$^*$Department of Mathematics, Northeastern University, Shenyang 110004, China\\
$^\dag$Department of Mathematical Sciences, University of Delaware, Delaware, 19716, USA}

\begin{abstract}
In this paper, we present and analyze a weak Galerkin finite element (WG) method for solving the symmetric hyperbolic systems. This method is highly flexible by allowing the use of discontinuous finite elements on element and its boundary independently of each other. By introducing special weak derivative, we construct a stable weak Galerkin scheme and derive the optimal $L_2$-error estimate of $O(h^{k+\frac{1}{2}})$-order for the discrete solution when the $k$-order polynomials are used for $k\geq 0$. As application, we discuss this WG method for solving the singularly perturbed convection-diffusion-reaction equation and derive an $\varepsilon$-uniform error estimate of order $k+1/2$. Numerical examples are provided to show the effectiveness of the proposed WG method.
\end{abstract}

\begin{keyword}
Weak Galerkin method; symmetric hyperbolic systems; stability; optimal error estimate; singularly perturbed problem

\MSC 65M60, 65N30, 65N12
\end{keyword}
\end{frontmatter}
\section{Introduction}
\setcounter{section}{1}\setcounter{equation}{0} In this work, we study a weak Galerkin finite element (WG) method for solving the first order hyperbolic systems:
\begin{eqnarray}
\sum_{k=1}^dA_k\partial_k\bo{u}+B{\bo u}={\bo f}, \;x\in \Omega\,,\la{1.1}
\end{eqnarray}
with the given boundary equation, see \rf{2.1}-\rf{2.2} for details. The corresponding scalar form of hyperbolic systems \rf{1.1} is the transport-reaction equation:
\begin{eqnarray}
\beta\cdot\nabla u+\alpha u=f,\;in\;\;\Omega,\;\;u=g,\;\;on\;\;\partial\Omega_-,\label{1.2}
\end{eqnarray}
where $\partial\Omega_-$ is the inflow boundary.

At present, for first order hyperbolic problems, one mainstream numerical method is the discontinuous Galerkin finite element (DG) method. It is well known that just for problem \rf{1.2}, the original DG method was presented by Reed and Hill in 1973 \cite{Reed} and was analyzed by
Lesaint and Raviart in \cite{Lesaint}. They showed that the DG method
has an explicit fashion if $\beta$ is constant and has the convergence rate
of $O(h^k)$-order when the method uses polynomials of order $k$. Later on,
Johnson and Pitkaranta \cite{Johnson} improved this convergence
order to $O(h^{k+\frac{1}{2}})$-order. Peterson in \cite{Peterson} further proved that the $
O(h^{k+\frac{1}{2}})$-order convergence is sharp within the quasi-uniform triangulation. However, under some special conditions on mesh and the convection direction $\beta$, Cockburn et. al. \cite{Cockburn,Cockburn1} and Richter \cite{Richter1} further showed that the $(k+1)$-order convergence rate can be obtained for the DG method. For symmetric hyperbolic systems \rf{1.1}, Ern and Guermond \cite{Ern,Ern1} made a systematic analysis for a class of DG methods. They gave the unique existence conditions of the exact solution of problem \rf{1.1} and established an abstract error estimate for the DG solution which implies that the convergence rate is of $O(h^{k+\frac{1}{2}})$-order if the $k$-order polynomials are used. Zhang et. al. in \cite{Zhang2,Zhang3} also proposed an upwind-like DG method for problem \rf{1.1}. Moreover, for the time-dependent problem of systems \rf{1.1}, Falk et. al. in \cite{Falk} and Monk et. al. in \cite{Monk} presented, respectively, the explicit and semi-explicit space-time DG methods.

As a new type of DG method, recently, the WG method has attracted much attention in the
field of numerical partial differential equations. This method was
introduced and analyzed originally by Wang and Ye in \cite{Wang} for second order elliptic
problems. Since then, WG methods have been used and analyzed for solving various partial differential equations, for example, the convection-diffusion-reaction equation, biharmonic equation, parabolic equation, Stokes equation and Navier-Stokes equation, and so on, see \cite{Xie,Gao,Lin1,Lin,Lijian,Lin2,Lin3,Wang,Wang2,Wang3,Zhang0,Zhang,Zhang1,Zhang4}. In general, a WG method can be considered as an extension of the standard
finite element method or DG method where classical derivatives are replaced in
the variational equation by the weakly defined derivatives on
discontinuous weak functions. The main feature of this method is: (1) the
weak finite element function $u_h=\{u_h^0,\,u_h^b\}$ is used in which $u_h^0$ is totally discontinuous on the partition and the value $u_h^b$ of $u_h$ on element edge may be
independent with the value $u_h^0$ of $u_h$ in the interior of element; (2)
the weak derivatives are introduced as distributions of weak finite element functions; The readers are referred to articles \cite{Lin1,Lin2,Wang2}
for more detailed explanation of this method and its relation with
other finite element methods. Although, WG methods have been studied for various partial differential equations, as authors' best knowledge, no WG method is presented for first order hyperbolic problems in existing literatures.

In this paper, we present and analyze a WG method imposed on shape regular meshes for solving symmetric hyperbolic systems \rf{1.1}. We first construct a stable WG scheme which has the feature that the unknown $\bo u_h^b$ in $\bo u_h=\{\bo u_h^0,\bo u_h^b\}$ can be eliminated locally, edge by edge, from the discrete WG equation so that the resulting in WG equation may be a linear system of equations only involving unknown ${\bo u}_h^0$. Therefore, the computation cost can be reduced greatly. Then, we do the error analysis for the WG method using the $k$-order polynomials for $k\geq 0$. We prove the following optimal error estimate for the WG solution $\bo u_h$:
\begin{equation}
|||\bo u-\bo u_h|||\leq Ch^{k+\frac{1}{2}}|\bo u|_{k+1},\;k\geq 0,\la{1.3}
\end{equation}
where $|||\cdot|||$ is an energy norm. As application of our method, besides the Maxwell's equations, we also consider the singularly perturbed convection-diffusion-reaction equation:
\begin{eqnarray}
-\varepsilon\triangle u+\boldsymbol{\beta}\cdot\nabla u+\alpha u=f,\;x\in\Omega,\;\;u=0,\;x\in \partial\Omega,\la{1.4}
\end{eqnarray}
where $\varepsilon>0$ is a parameter. It is well known that when the diffusion coefficient $\varepsilon$ is very small, the solution of this boundary value problem typically possesses $layers$, which is thin regions where the solution and/or its derivatives change rapidly. Standard numerical methods fail to provide accurate approximations in this case unless the computational mesh is of the magnitude of the layers.

We use the proposed WG method to solve problem \rf{1.4} by transforming it into the symmetric hyperbolic systems \rf{1.1} and obtain the error estimate:
\begin{equation}
\|u-u_h^0\|+\sqrt{\varepsilon}\,\|\nabla u-\nabla^wu_h\|\leq Ch^{k+\frac{1}{2}}(|u|_{k+1}+\sqrt{\varepsilon}\,|\nabla u|_{k+1}),\,k\geq 0,\la{1.5}
\end{equation}
where $\nabla^wu_h$ is the weak gradient approximation of $\nabla u$ defined in an appropriate manner. Our error estimate holds true uniformly with respect to $\varepsilon$ and without imposing extra conditions on the mesh and the data $\beta$ and $\alpha$, see \cite{Xie}. On the other hand, Lin et. al. in \cite{Lin} also presented a WG method for solving problem \rf{1.4} directly and derived the following error estimate:
\begin{equation}
\|u-u_h^0\|+\sqrt{\varepsilon}\,\|\nabla u-\nabla_wu_h\|\leq Ch^k(h^{\frac{1}{2}}+\sqrt{\varepsilon}\,)\,|u|_{k+1},\,k\geq 0,\la{1.6}
\end{equation}
where $\nabla_wu_h$ is the weak gradient \cite{Wang}. Obviously, our result \rf{1.5} is better than result \rf{1.6}. For example, for the piecewise constant element ($k=0$),
our WG method gives an $O(h^{\frac{1}{2}})$-order convergence rate, but no convergence rate can be obtained from estimate \rf{1.6} if $k=0$.

Our work provides an approach to develop the WG method for first order hyperbolic problems.

The rest of this paper is organized as follows. In Section 2, we introduce the symmetric hyperbolic systems and construct the corresponding WG scheme. In Section 3, we show the stability of this WG scheme and derive the optimal $L_2$-error estimate of $O(h^{k+\frac{1}{2}})$-order for $k\geq 0$. To show the effectiveness of the proposed WG method,  numerical experiments are provided in Section 4 for solving the Maxwell's equations and the singularly perturbed convection-diffusion-reaction equation.

Throughout this paper, we adopt the notations $H^m(D)$ to
indicate the usual Sobolev spaces on subdomain $D\subset \Omega$
equipped with the norm $\|\cdot\|_{m,D}$ and semi-norm
$|\cdot|_{m,D}$, and when $D=\Omega$, we omit the
index $D$. The inner product and norm in space $H^0(\Omega)=L_2(\Omega)$ are
denoted by $(\cdot,\cdot)$ and $\|\cdot\|$, respectively. We
use letter $C$ to represent a generic powsitive constant,
independent of the mesh size $h$.
\section{Symmetric hyperbolic systems and its weak Galerkin approximation}
\setcounter{section}{2}\setcounter{equation}{0}
 Consider the following first-order hyperbolic
system:
\begin{eqnarray}
\sum_{k=1}^dA_k\partial_k\bo{u}+B{\bo u}={\bo f}, \;x\in \Omega\,,\la{2.1}\\
(M-D_n){\bo u}={\bo 0},\;x\in
\partial\Omega.\la{2.2}
\end{eqnarray}
Here, $\Omega\subset R^d$ is a bounded polygonal or polyhedral domain, $A_k=(a_{ij}^{(k)}(x))$,
$k=1,\cdots,d,\,B=(b_{ij}(x))$ and $M=(m_{ij}(x))$ are some given
$m\times m$ matrices, $D_n=\sum_{k=1}^dA_kn_k$, $n(x)=(n_1,\cdots,n_d)^T$ is the outward unit normal vector
at the point $x\in
\partial\Omega$, ${\bo u}=(u_1,\cdots,u_m)^T$ and ${\bo
f}=(f_1,\cdots,f_m)^T$ are
$m$-dimensional vector functions. In what follows, for simplicity, we denote by $\boldsymbol{A}=(A_1,\cdots,A_d)^T$ the vector matrix function and set
$$
\bo{A}\cdot\nabla{\bo u}=\sum_{k=1}^dA_k\partial_k{\bf u},\;\;\hbox{div}\bo{A}=\partial_1A_1+\cdots+\partial_d A_d,\;D_n={\bo A}\cdot n=\sum_{k=1}^d=A_kn_k.
$$
We assume that problem
\rf{2.1}-\rf{2.2} is a positive and symmetric hyperbolic system,
namely,
\begin{eqnarray}
&&A_i=A_i^T,\;i=1,\cdots,d,\;x\in\Omega,\label{2.3}\\
&&B+B^T-\hbox{div}\boldsymbol{A}\geq 2\sigma_0I,\; x\in \Omega,\label{2.4}\\
&&M+M^T\geq 0,\;x\in \partial\Omega,\label{2.5}\\
&&ker(M-D_n)+ker(M+D_n)=R^m,\;x\in \partial\Omega,\label{2.6}
\end{eqnarray}
where constant $\sigma_0>0$ and by using the
expression $A\geq 0$ we imply that matrix $A$ is
positive semi-definite. Under the assumptions
\rf{2.3}-\rf{2.6}, according to Friedrichs' theory \cite{Friedrichs}, problem \rf{2.1}-\rf{2.2} has a
unique solution under appropriate smoothness conditions, also see \cite{Ern,Ern1}. It should be pointed that condition \rf{2.6} is not used in our stability and error analysis below.

Problem \rf{2.1}-\rf{2.2} can describe many important physics
processes. Two examples of such symmetric hyperbolic systems are as follows.

{\bf Maxwell's equations}

Let $\nu$ and $\sigma$
be two positive functions in $L_\infty(\Omega)$ uniformly bounded
away from zero. Consider the following Maxwell's equations in
$R^3$
\begin{eqnarray}
\nu H+\nabla\times E=h,\;x\in\Omega,\la{2.7}\\
\sigma E-\nabla\times H=g,\;x\in\Omega,\la{2.8}\\
 E\times n=0,\;x\in\partial\Omega,\la{2.9}
\end{eqnarray}
where $H$, $E$, $h$ and $g$ are three-dimensional vector functions. This
problem can be cast into the form of a positive and symmetric
hyperbolic system by setting ${\bf u}=(H,E)^T$,
$$
A_k=\left(
\begin{array}{cc}
   O&Q_k \\
  Q_k^T &O \\
\end{array}
\right),_{_{\; \displaystyle{k=1,2,3,}}}\;\; B=\left(
\begin{array}{cc}
  \nu I &O \\
  O  & \sigma I \\
\end{array}
\right),\;\; {\bf f}=\left(
\begin{array}{c}
  h\\
  g\\
\end{array}
\right),
$$
$$
Q_1=\left(
\begin{array}{ccc}
  0 & 0 & 0 \\
  0 & 0 & -1 \\
  0 &1 & 0 \\
\end{array}
\right),\;\; Q_2=\left(
\begin{array}{ccc}
  0 & 0 & 1 \\
  0 & 0 & 0 \\
  -1 & 0 & 0 \\
\end{array}
\right),\;\; Q_3=\left(
\begin{array}{ccc}
  0 & -1 & 0 \\
  1 & 0 & 0 \\
  0 & 0 & 0 \\
\end{array}
\right),
$$
and choosing the boundary matrix
$$
M=\left(
\begin{array}{cc}
  O & -R \\
  R^T & O \\
\end{array}
\right),\;\; R=\sum_{k=1}^3Q_kn_k,
$$
where $O$ is the $3\times 3$ zero matrix. The conditions
\rf{2.3}-\rf{2.6} can be verified directly in which $\sigma_0=\min\{\nu,\sigma\}$.

{\bf Convection-diffusion-reaction problem}

Consider the singularly perturbed convection-diffusion-reaction problem in $R^2$:
\begin{eqnarray}
-\varepsilon\triangle u+\boldsymbol{\beta}\cdot\nabla u+\alpha u=f,\;x\in\Omega,\la{2.10}\\
u=0,\;x\in \partial\Omega,\la{2.11}
\end{eqnarray}
where $\varepsilon>0$ is a small parameter and $\alpha-\hbox{div}\beta/2\geq\alpha_0>0$. This problem can be written as a first order hyperbolic system:
\begin{eqnarray}
\boldsymbol{\sigma}+\sqrt{\varepsilon}\,\nabla u=\boldsymbol{0},\;x\in \Omega,\la{2.11a}\\
\sqrt{\varepsilon}\,\hbox{div}\boldsymbol{\sigma}+\boldsymbol{\beta}\cdot\nabla u+\alpha u=f,\;x\in\Omega,\la{2.11b}\\
u=0,\;x\in\partial\Omega,\la{2.11c}
\end{eqnarray}
which has the form of positive and symmetric, hyperbolic system by
setting ${\bo u}=(\sigma_1,\sigma_2,u)^T$,
$$
A_1=\left(
\begin{array}{ccc}
  0 & 0 & \sqrt{\varepsilon} \\
  0 & 0 & 0 \\
  \sqrt{\varepsilon} & 0 & \beta_1 \\
\end{array}
\right),\;\; A_2=\left(
\begin{array}{ccc}
  0 & 0 & 0 \\
  0 & 0 & \sqrt{\varepsilon} \\
  0 & \sqrt{\varepsilon} & \beta_2 \\
\end{array}
\right),\;\; B=\left(
\begin{array}{ccc}
  1 & 0 & 0 \\
  0 & 1 & 0 \\
  0 & 0 & \alpha \\
\end{array}
\right),\;\; {\bo f}=\left(
\begin{array}{c}
  0 \\
  0 \\
  f \\
\end{array}
\right),
$$
and choosing the boundary matrix
$$
M=\left(
\begin{array}{ccc}
  0&0& -\sqrt{\varepsilon}n_1 \\
  0&0&-\sqrt{\varepsilon}n_2 \\
  \sqrt{\varepsilon}n_1&\sqrt{\varepsilon}n_2&1\\
\end{array}
\right),\;\;
D_n=\left(
\begin{array}{ccc}
  0&0& \sqrt{\varepsilon}n_1 \\
  0&0&\sqrt{\varepsilon}n_2 \\
  \sqrt{\varepsilon}n_1&\sqrt{\varepsilon}n_2&\beta\cdot n\\
\end{array}
\right).
$$
The conditions \rf{2.3}-\rf{2.6} can be verified directly in which $\sigma_0=\min\{1,\alpha_0\}$.

In the above examples, although the boundary matrices $\{M\}$
should be determined by the boundary value conditions of the
problems, they are not unique.
\commentout{Here we have chosen the boundary
matrices $\{M\}$ properly such that they also satisfy our
argument requirement in the error analysis, see \rf{3.16}.
}

Now we introduce the WG method for problem \rf{2.1}-\rf{2.2}.

Let $T_h=\bigcup\lbrace K \rbrace$ be a partition of
domain $\Omega$ that consists of arbitrary polygons/polyhedra, where the mesh size $h=\max \,h_K$, $h_K$ is the diameter of element $K$. Assume that the partition $T_h$ is shape regular defined by a set of conditions given in \cite{Wang2}.

First, let us recall the concepts of weak function and weak finite element space (see, e.g.,\cite{Wang,Wang2}) which will then be employed to define a weak Galerlin finite element scheme for problem \rf{2.1}-\rf{2.2}.
A weak function on element $K$ refers to a function $v=\{ v^0,v^b\}$ with $v^0=v|_{K}\in
L_2(K)$ and $v^b=v|_{\partial K}\in L_2(\partial K)$. Note that for
a weak function $v=\{v^0,v^b\}$, $v^b$ may not be necessarily the trace of $v^0$ on element boundary $\partial K$.

Introduce the weak Galerkin finite element spaces on partition $T_h$:
\begin{eqnarray*}
V_h=\{v=\{v^0,v^b\}:\;v|_K\in P_k(K),\,v^b|_e\in P_k(e),\,e\subset\partial K,\; K\in
T_h\},\;k\geq 0,
\end{eqnarray*}
where $P_k(D)$ is the space composed of all polynomials on a set $D$ with degree no more than
$k$. We emphasize that, for $v=\{v^0,v^b\}\in V_h$, $v^b$ is single valued on edge/face $e\subset\partial K$ which means that $v^b$ is continuous across $\partial K$. On the other hand, the component $v^0$ is defined element-wise and completely discontinuous on $T_h$. In a certain sense, a weak finite element function $v=\{v^0,v^b\}\in V_h$ is
formed with its components inside all elements glued together by its components on all edges/faces.

Different from the weak gradient and weak divergence usually used in the WG method, we introduce here a special weak derivative for our problem.
\commentout{
In order to easily understand our definition, we first consider the weak directional  derivative for a scalar function $u$. Let $\beta=(\beta_1,\cdots,\beta_d)^T$ be a direction vector and $\nabla_\beta u=\beta\cdot\nabla u$ be the directional derivative. By the Green's formula, we have
$$
(\nabla_\beta u,q)_K=-(u,\beta\cdot\nabla q+\hbox{div}\beta\, q)_K+\int_{\partial K}\beta\cdot nu\,qds,\;q\in H^1(K).
$$
Then, for given weak function $v=\{v^0,v^b\}\in V_h$, we can define its weak directional derivative $\nabla_{w,\beta}v|_K\in P_k(K)$ as the unique solution of the following equation:
\begin{equation}
(\nabla_{w,\beta} v,q)_K=-(v^0,\beta\cdot\nabla q+\hbox{div}\beta q)_K+\int_{\partial K}\beta\cdot{ n}v^bqds,\;\forall\,q\in P_k(K),\;K\in T_h.\la{2.12}
\end{equation}
}
For vector function ${\bo u}\in [H^1(K)]^m$ and vector matrix ${\bo A}=(A_1,\cdots,A_d)^T$, by using the Green's formula,  we have
$$
({\bo A}\cdot\nabla \bo{u,q})_K=-({\bo u}^0,{\bo A}\cdot\nabla {\bo q}+\hbox{div}{\bo A}\,{\bo q})_K+\int_{\partial K}D_n{\bo  u}\cdot{\bo q}ds,\;\forall\,{\bo q}\in [H^1(K)]^m,
$$
where $D_n={\bo A}\cdot n$. Therefore, for weak vector function ${\bo v}=\{{\bo v}^0,{\bo v}^b\}\in [V_h]^m$, we define its weak derivative $\nabla_{w,A}{\bo v}|_K\in [P(K)]^m$ related to vector matrix ${\bo A}$ on element $K\in T_h$ as the unique solution of the following equation:
\begin{equation}
(\nabla_{w,A}\bo{v,q})_K=-({\bo v}^0,{\bo A}\cdot\nabla {\bo q}+\hbox{div}{\bo A}\,{\bo q})_K+\int_{\partial K}D_n{\bo  v}^b\cdot{\bo q}ds,\;\forall\,{\bo q}\in [P_k(K)]^m.\la{2.13}
\end{equation}
Obviously, operator $\nabla_{w,A}$ is an analogy of the differential operator ${\bo A}\cdot \nabla$.

Denote the set $\partial T_h=\{\partial K:\,K\in T_h\}$. For simplicity, we use the following notations,
\begin{eqnarray*}
&&(w,v)_{T_h}=\sum_{K\in T_h}(w,v)_K,\;\;\;\langle w,v\rangle_{\partial T_h}=\sum_{K\in T_h}\langle w,v\rangle_{\partial K}=\sum_{K\in T_h}\int_{\partial K}wvds,\\
&&\langle w,v\rangle_{\partial \Omega}=\int_{\partial\Omega}uvds=\sum_{K\in T_h}\int_{\partial\Omega\cap\partial K}wvds.
\end{eqnarray*}
Let ${\bo u}\in [H^1(\Omega)]^m$ be the exact solution of problem \rf{2.1}-\rf{2.2}. Then, ${\bo u}$ satisfies the variational equation, for ${\bo q}\in [L_2(\Omega)]^m,\,\widehat{\bo q}\in [L_2(\partial\Omega)]^m$,
\begin{equation}
({\bo A}\cdot\nabla{\bo u},{\bo q})+(B{\bo u},{\bo q})+\frac{1}{2}\langle (M-D_n)\bo u, \widehat{\bo q}\rangle_{\partial\Omega}=({\bo f,\bo q}).\la{2.14}
\end{equation}
Motivated by this weak form, we introduce the bilinear form for ${\bo w}=\{{\bo w}^0,{\bo w}^b\},{\bo v}=\{{\bo v}^0,{\bo v}^b\}\in [V_h]^m$,
\begin{equation}
a({\bo w,\bo v})=(\nabla_{w,A}{\bo w,\bo v}^0)_{T_h}+(B{\bo w}^0,{\bo v}^0)_{T_h}+\frac{1}{2}\langle (M-D_n){\bo w}^b,{\bo v}^b\rangle_{\partial\Omega}.\la{2.15}
\end{equation}
{\sc Weak Galerkin Method}: a weak Galerkin finite element approximation for problem \rf{2.1}-\rf{2.2} is to find ${\bo u}_h=\{{\bo u}_h^0,{\bo u}_h^b\}\in [V_h]^m$ such that
\begin{equation}
a({\bo u}_h,{\bo v})+s({\bo u}_h,{\bo v})=({\bo f},{\bo v}^0)_{T_h},\;\forall\,{\bo v}=\{{\bo v}^0,{\bo v}^b\}\in [V_h]^m,\la{2.16}
\end{equation}
where the stabilizer
$$
s({\bo w,\bo v})=\langle \mu({\bo w}^0-{\bo w}^b),{\bo v}^0-{\bo v}^b\rangle_{\partial T_h},
$$
and the parameter $\mu>0$ can be chosen properly to enhance the stability of this WG scheme.
In what follows, we always choose $\mu$ such that
\begin{equation}
\mu-\frac{1}{2}\rho(D_n(x))\geq\mu_0>0,\;x\in \overline{\Omega},\la{2.17}
\end{equation}
where $\rho(D_n)=\|D_n\|_2$ represents the spectral radius of matrix $D_n$.

Below let us give a discussion on the solving method of WG equation \rf{2.16}. Since WG equations concern the unknown ${\bo u}^0_h$ on elements and ${\bo u}^b_h$ on element boundaries, it seems that the number of unknowns of a WG equation  is much more than that of the usual finite element equation. But, WG equations usually have a hybridized construction so that the unknown ${\bo u}^b_h$ can be eliminated locally by means of the unknown ${\bo u}_h^0$. Therefore, the WG equation can yield a system of equations involving much less number of unknowns than what it appears. To show this, let us examine the WG equation \rf{2.16}. By using weak derivative formula \rf{2.13}, we can write equation \rf{2.16} in the following form
\begin{eqnarray*}
&&-({\bo u}_h^0,{\bo A}\cdot\nabla{\bo v}^0)_{T_h}+((B-\hbox{div}{\bo A}){\bo u}_h^0,{\bo v}^0)_{T_h}+\langle D_n{\bo u}^b_h,{\bo v}^0\rangle_{\partial T_h}\nonumber\\
&&+\frac{1}{2}\langle (M-D_n){\bo u}_h^b,{\bo v}^b\rangle_{\partial\Omega}+\langle \mu({\bo u}_h^0-{\bo u}_h^b),{\bo v}^0-{\bo v}^b\rangle_{\partial T_h}=
({\bo f},{\bo v}^0),\;\forall\,{\bo v}\in [V_h]^m.
\end{eqnarray*}
Taking ${\bo v}=\{{\bo v}^0,0\}$ and ${\bo v}=\{0,{\bo v}^b\}$, respectively, it yields
\begin{eqnarray}
&-&({\bo u}_h^0,{\bo A}\cdot\nabla{\bo v}^0)_{T_h}+((B-\hbox{div}{\bo A}){\bo u}_h^0,{\bo v}^0)_{T_h}
+\langle D_n{\bo u}^b_h,{\bo v}^0\rangle_{\partial T_h}\nonumber\\
&&+\langle \mu({\bo u}_h^0-{\bo u}_h^b),{\bo v}^0\rangle_{\partial T_h}=({\bo f},{\bo v}^0),\; \la{2.18}\\
&&\frac{1}{2}\langle (M-D_n){\bo u}_h^b,{\bo v}^b\rangle_{\partial\Omega}-\langle \mu({\bo u}_h^0-{\bo u}_h^b),{\bo v}^b\rangle_{\partial T_h}=
0.\;\;\la{2.19}
\end{eqnarray}
Now, for any fixed edge/face $e\subset\partial K,\,K\in T_h$, taking ${\bo v}^b|_e\in [P_k(e)]^m$ and ${\bo v}^b=0$ on  other edges/faces, we obtain from equation \rf{2.19} that
\begin{eqnarray}
\left\{\begin{array}{ll}
\displaystyle{\frac{1}{2}}\langle (M-D_n){\bo u}_h^b,{\bo v}^b\rangle_{e}+\langle \mu{\bo u}_h^b,{\bo v}^b\rangle_{e}=\langle \mu{\bo u}_h^0,{\bo v}^b\rangle_{e},\,e\subset\partial\Omega,\\
2\langle {\bo u}_h^b,{\bo v}^b\rangle_{e}=\langle {\bo u}_h^0,{\bo v}^b\rangle_{e\cap\partial K}+\langle {\bo u}_h^0,{\bo v}^b\rangle_{e\cap\partial K'},\,e\not\subset\partial\Omega,\,{\bo v}^b\in [P_k(e)]^m,\la{2.20}
\end{array}
 \right.
\end{eqnarray}
where $K$ and $K'$ are two adjacent elements sharing the common edge $e$. By conditions \rf{2.5} and \rf{2.17}, we have for any $\bo q\in [L_2(e)]^m$,
$$
\langle M{\bo q},{\bo q}\rangle_e=\frac{1}{2}\langle (M+M^T){\bo q},{\bo q}\rangle_e\geq 0,\;\langle (\mu I-\frac{1}{2}D_n){\bo q},{\bo q}\rangle_e\geq \mu_0\langle{\bo q},{\bo q}\rangle_e,\,e\subset\partial\Omega.
$$
Hence, for each edge/face $e\subset\partial K,\,K\in T_h$, ${\bo u}_h^b|_e$ can be solved uniquely from equation \rf{2.20} by means of ${\bo u}_h^0(K)$ and ${\bo u}_h^0(K')$. Thus, the unknown $\bo u^b_h$ can be eliminated from equation \rf{2.18} so that WG equation \rf{2.18} may be a linear system of equations only involving unknown ${\bo u}_h^0$. After solving $\bo u_h^0$, we can solve $\bo u_h^b$ from equation \rf{2.20}, edge by edge. Actually, in application, we usually only need to find solution $\bo u_h^0$. Therefore, in general, a WG method has a comparable computation cost with the DG method.

\section{Stability and error analysis}
\setcounter{section}{3}\setcounter{equation}{0}
In this section, we establish the stability of WG scheme \rf{2.16} and give the error estimate for the WG solution.
\subsection{Stability}
We first give a lemma.\\
\yl{ 3.1}\quad{\em For ${\bo v}=\{ {\bo v}^0,{\bo v}^b\}\in [V_h]^m$, it holds
\begin{eqnarray}
a({\bo v,\bo v})+s({\bo v,\bo v})&=&\frac{1}{2}(N{\bo v}^0,{\bo v}^0)_{T_h}+\langle(\mu I-\frac{1}{2}D_n)({\bo v}^0-{\bo v}^b),{\bo v}^0-{\bo v}^b\rangle_{\partial T_h}\nonumber\\
&&+\frac{1}{2}\langle M{\bo v}^b,{\bo v}^b\rangle_{\partial\Omega}.\la{3.1}
\end{eqnarray}
where matrix} $N=B+B^T-\hbox{div}{\bo A}$.\\
\zm\quad By using weak derivative formula \rf{2.13} and the Green's formula, we have for ${\bo w,\bo v}\in [V_h]^m$,
\begin{eqnarray*}
&&(\nabla_{w,A}{\bo w},\bo{v}^0)_K=-({\bo w}^0,{\bo A}\cdot\nabla{\bo v}^0)_{K}-(\hbox{div}{\bo A}{\bo w}^0,{\bo v}^0)_{K}+\langle D_n{\bo w}^b,{\bo v}^0\rangle_{\partial K}\\
&=&({\bo A}\cdot\nabla{\bo w}^0,{\bo v}^0)_{K}+(\hbox{div}{\bo A}{\bo w}^0,{\bo v}^0)_{K}-\langle D_n{\bo v}^0,{\bo w}^0\rangle_{\partial K}\\
&&-(\hbox{div}{\bo A}{\bo w}^0,{\bo v}^0)_{K}+\langle D_n{\bo w}^b,{\bo v}^0\rangle_{\partial K}.
\end{eqnarray*}
Again using formula \rf{2.13}, it yields
$$
({\bo A}\cdot\nabla{\bo w}^0,{\bo v}^0)_{K}+(\hbox{div}{\bo A}{\bo w}^0,{\bo v}^0)_{K}=-(\nabla_{w,A}{\bo v},\bo{w}^0)_K
+\langle D_n{\bo v}^b,{\bo w}^0\rangle_{\partial K}.
$$
Combining the above two equalities, we obtain (noting that $D_n$ is symmetric matrix)
\begin{eqnarray*}
(\nabla_{w,A}{\bo w},\bo{v}^0)_K&=&-(\nabla_{w,A}{\bo v},\bo{w}^0)_K
+\langle D_n({\bo v}^b-{\bo v}^0),{\bo w}^0\rangle_{\partial K}-(\hbox{div}{\bo A}{\bo w}^0,{\bo v}^0)_{K}\\
&&+\langle D_n{\bo w}^b,{\bo v}^0\rangle_{\partial K}\\
&=&-(\nabla_{w,A}{\bo v},\bo{w}^0)_K-(\hbox{div}{\bo A}{\bo w}^0,{\bo v}^0)_{K}
+\langle D_n({\bo v}^b-{\bo v}^0),{\bo w}^0-{\bo w}^b\rangle_{\partial K}\\
&&+\langle D_n({\bo v}^b-{\bo v}^0),{\bo w}^b\rangle_{\partial K}
+\langle D_n{\bo w}^b,{\bo v}^0\rangle_{\partial K}\\
&=&-(\nabla_{w,A}{\bo v},\bo{w}^0)_K-(\hbox{div}{\bo A}{\bo w}^0,{\bo v}^0)_{K}
-\langle D_n({\bo v}^0-{\bo v}^b),{\bo w}^0-{\bo w}^b\rangle_{\partial K}\\
&&+\langle D_n{\bo v}^b,{\bo w}^b\rangle_{\partial K}.
\end{eqnarray*}
Taking $\bo w=\bo v$ and summing for $K\in T_h$, it yields
\begin{eqnarray*}
(\nabla_{w,A}{\bo v},\bo{v}^0)_{T_h}=-\frac{1}{2}(\hbox{div}{\bo A}{\bo v}^0,{\bo v}^0)_{T_h}
-\frac{1}{2}\langle D_n({\bo v}^0-{\bo v}^b),{\bo v}^0-{\bo v}^b\rangle_{\partial T_h}
+\frac{1}{2}\langle D_n{\bo v}^b,{\bo v}^b\rangle_{\partial\Omega}.
\end{eqnarray*}
where we have used the fact that since ${\bo v}^b$ is continuous across $\partial K$ and on common edge $e=\partial K\cap\partial K'$, $D_n=-D_{n'}$, so that $\langle D_n{\bo v}^b,{\bo v}^b\rangle_{\partial T_h}=\langle D_n{\bo v}^b,{\bo v}^b\rangle_{\partial\Omega}$ holds. Hence, it follows from \rf{2.15} that
\begin{eqnarray}
&&a({\bo v,\bo v})+s({\bo v,\bo v})=\frac{1}{2}((B+B^T-\hbox{div}{\bo A}){\bo v}^0,{\bo v}^0)_{T_h}\nonumber\\
&+&\langle(\mu I-\frac{1}{2}D_n)({\bo v}^0-{\bo v}^b),{\bo v}^0-{\bo v}^b\rangle_{\partial T_h}+\frac{1}{2}\langle M{\bo v}^b,{\bo v}^b\rangle_{\partial\Omega}.\la{3.2}
\end{eqnarray}
The proof is completed.\zb

Introduce the notation:
\begin{eqnarray}
|||\bo v|||^2=\sigma_0(\bo{v}^0,\bo{v}^0)_{T_h}+\mu_0\langle\bo{v}^0-\bo{v}^b,\bo{v}^0-\bo{v}^b\rangle_{\partial T_h}+\frac{1}{2}\langle M\bo{v}^b,\bo{v}^b\rangle_{\partial\Omega}.\;\la{3.3}
\end{eqnarray}
It is easy to see that $|||{\bo v}|||$ defines a norm on space $[V_h]^m$. In fact, when $|||{\bo v}|||=0$, we obtain ${\bo v}^0=0$ and $\langle {\bo v}^b,{\bo v}^b\rangle_{\partial T_h}=0$, so $\bo{v}=\{\bo{v}^0,\bo{v}^b\}=0$ holds.\\
\dl{ 3.1}\quad{\em The weak Galerkin finite element equation \rf{2.16} has one unique solution $\bo{u}_h\in [V_h]^m$ and the following stability estimate holds.}
\begin{equation}
|||{\bo u}_h|||\leq \frac{1}{\sqrt{\sigma_0}}\|\bo f\|.\la{3.4}
\end{equation}
\zm\quad We only need to prove the stability estimate \rf{3.4}. By Lemma 3.1, \rf{2.4}-\rf{2.5} and \rf{2.17}, we first obtain
\begin{equation}
|||{\bo v}|||^2\leq a(\bo v,\bo v)+s(\bo v,\bo v),\;\forall\, {\bo v}\in [V_h]^m.\la{3.4a}
\end{equation}
Then, taking ${\bo v}={\bo u_h}$ in equation \rf{2.16} and using inequality: $\sqrt{\sigma_0}\,\|{\bo u}_h^0\|\leq |||\bo {u}_h|||$, estimate \rf{3.4} is derived.\zb
\subsection{Error analysis}

Given a function $u$ with sufficient regularity, one can find an approximation to $u$ by either interpolation or projection in a standard finite element space. In the WG space $V_h$, we will use a locally defined projection of $u$ as its basic approximation. Specifically, let $P_h:u\in L_2(K)\rightarrow P_h u\in P_k(K)$
be the local $L_2$ projection operator such that
\begin{equation}
(u-P_hu,q)_{K}=0,\;\forall\,q\in P_k(K),\,K\in T_h.\la{3.5}
\end{equation}
Operator $P_h$ has the approximation property:
\begin{equation}
\|u-P_hu\|_{s,K}\leq Ch_K^{m-s}\|u\|_{m,K},\;0\leq s\leq m\leq k+1.\la{3.6}
\end{equation}
The $L_2$ projection operator $P_{\partial K}: L_2(e)\rightarrow P_{k}(e),\,e\subset \partial K$ can be defined similarly on the edges of element
$K\in T_h$. Now, we define a projection operator $Q_h: u\in H^1(\Omega)\rightarrow Q_hu\in V_h$ by its action on each element $K$ such that
\begin{equation}
Q_hu|_K=\{Q^0u,Q^bu\}\doteq\{P_hu,P_{\partial K}u\},\;K\in T_h.\la{3.7}
\end{equation}
For vector function ${\bo u}\in [H^1(\Omega)]^m$, we set $Q_h{\bo u}=(Q_hu_1,\cdots,Q_hu_m)^T$. \\
\yl{ 3.2}\quad{\em Let ${\bo u}\in [H^1(\Omega)]^m$. Then it holds true for $\bo v=\{\bo{v}^0,\bo{v}^b\}\in [V_h]^m$,
\begin{equation}
({\bo A}\cdot\nabla \bo{u},{\bo v}^0)_{T_h}=(\nabla_{w,A}Q_h{\bo u},{\bo v}^0)_{T_h}+\langle D_n(\bo u-Q^b\bo u),\bo v^0\rangle_{\partial T_h}+l_1({\bo u},\bo{v}),\la{3.8}
\end{equation}
where}
\begin{eqnarray}
l_1(\bo u,\bo v)=(Q^0\bo u-\bo u,{\bo A}\cdot\nabla{\bo v}^0)_{T_h}+(\hbox{div}\bo A(Q^0\bo u-\bo u),\bo v^0)_{T_h}.\la{3.9}
\end{eqnarray}
\zm\quad From the Green's formula and weak derivative formula \rf{2.13}, we have
\begin{eqnarray*}
({\bo A}\cdot\nabla{\bo u},v^0)_{T_h}&=&-(\bo u,\bo A\cdot\nabla{\bo v}^0)_{T_h}-(\hbox{div}\bo A\bo u,\bo v^0)_{T_h}+
\langle D_n\bo u,\bo v^0\rangle_{\partial T_h}\\
&=&-(Q^0\bo u,\bo A\cdot\nabla{\bo v}^0)_{T_h}-(\hbox{div}\bo AQ^0\bo u,\bo v^0)_{T_h}+l_1(\bo u,\bo v)+
\langle D_n\bo u,\bo v^0\rangle_{\partial T_h}\\
&=&(\nabla_{w,A}Q_h\bo u,{\bo v}^0)_{T_h}-\langle D_nQ^b\bo u,\bo v^0\rangle_{\partial T_h}+l_1(\bo u,\bo v)+
\langle D_n\bo u,\bo v^0\rangle_{\partial T_h}.
\end{eqnarray*}
The proof is completed.\zb

Set
\begin{eqnarray}
&&l_2(\bo u,\bo v)=(B(\bo u-Q^0\bo u),\bo v^0)_{T_h},\la{3.10}\\
&&l_3(\bo u,\bo v)=\langle D_n(\bo u-Q^b\bo u),\bo v^0-\bo v^b\rangle_{\partial T_h},\la{3.11}\\
&&l_4(\bo u\bo ,v)=\frac{1}{2}\langle(M+D_n)(\bo u-Q^b\bo u),\bo v^b\rangle_{\partial\Omega}.\la{3.12}
\end{eqnarray}
\yl{ 3.3}\quad{\em Let ${\bo u}\in [H^1(\Omega)]^m$ be the solution of problem \rf{2.1}-\rf{2.2}. Then we have for $\bo v=\{\bo v^0,\bo v^b\}\in [V_h]^m$},
\begin{eqnarray}
a(Q_h\bo u,\bo v)=(\bo f,\bo v^0)-l_1(\bo u,\bo v)-l_2(\bo u,\bo v)-l_3(\bo u,\bo v)-l_4(\bo u,\bo v).\la{3.13}
\end{eqnarray}
\zm\quad From equations \rf{2.1}-\rf{2.2}, we have for $\bo v=\{\bo v^0,\bo v^b\}\in [V_h]^m$,
\begin{eqnarray*}
(\bo A\cdot\nabla{\bo u},{\bo v}^0)_{T_h}+(B{\bo u},{\bo v^0})+\frac{1}{2}\langle (M-D_n){\bo u,\bo v^b}\rangle_{\partial\Omega}=({\bo f,\bo v^0}),
\end{eqnarray*}
together with Lemma 3.2, it yields
\begin{eqnarray}
&&(\nabla_{w,A}Q_h{\bo u},{\bo v}^0)_{T_h}+\langle D_n(\bo u-Q^b\bo u),\bo v^0\rangle_{\partial T_h}+l_1({\bo u},\bo{v})\nonumber\\
&+&(B{\bo u},{\bo v^0})+\frac{1}{2}\langle (M-D_n){\bo u,\bo v^b}\rangle_{\partial\Omega}=({\bo f,\bo v^0}).\la{3.14}
\end{eqnarray}
Since
\begin{eqnarray*}
&&(B{\bo u},{\bo v^0})=(BQ^0{\bo u},{\bo v^0})+l_2(\bo u,\bo v),\\
&&\frac{1}{2}\langle (M-D_n){\bo u,\bo v^b}\rangle_{\partial\Omega}=\frac{1}{2}\langle (M-D_n)Q^b\bo u,\bo v^b\rangle_{\partial\Omega}+\frac{1}{2}\langle (M-D_n)(\bo u-Q^b\bo u),\bo v^b\rangle_{\partial\Omega},
\end{eqnarray*}
then we have from \rf{3.14} and the definition \rf{2.15} of $a(\bo w,\bo v)$ that
\begin{eqnarray}
&&a(Q_h{\bo u},{\bo v})+l_2({\bo u},\bo{v})+\langle D_n(\bo u-Q^b\bo u),\bo v^0\rangle_{\partial T_h}+l_1({\bo u},\bo{v})\nonumber\\
&+&\frac{1}{2}\langle (M-D_n)(\bo u-Q^b\bo u),\bo v^b\rangle_{\partial\Omega}=({\bo f,\bo v^0}).\la{3.15}
\end{eqnarray}
Noting that $\bo u-Q^b\bo u$ and $\bo v^b$ are continuous across element boundaries, it implies
$$
\langle D_n(\bo u-Q^b\bo u),\bo v^0\rangle_{\partial T_h}=\langle D_n(\bo u-Q^b\bo u),\bo v^0-\bo v^b\rangle_{\partial T_h}+\langle D_n(\bo u-Q^b\bo u),\bo v^b\rangle_{\partial\Omega}.
$$
Substituting this into \rf{3.15}, the proof is completed.\zb
\commentout{
In order to do the error analysis, we still need an assumption on problem \rf{2.1}-\rf{2.2}.
\begin{eqnarray}
(A)\quad &&\hbox{On each edge/face $e\subset\partial\Omega$, either matrices $M$ and $D_n$ are constant or there }\nonumber\\
&&\hbox{exists a constant $C_M>0$ such that}\nonumber\\
&&|\langle(M+D_n)\bo u,\bo w\rangle_{\partial\Omega}|\leq C_M\|\bo u\|_{L_2(\partial\Omega}\langle M\bo w,\bo w\rangle^{\frac{1}{2}}_{\partial\Omega},\;\forall\,\bo u,\bo w\in [L_2(\partial\Omega)]^m.\la{3.16}
\end{eqnarray}
Assumption (A) can be satisfied by many practical problems. For example, for the Maxwell's equations \rf{2.7}-\rf{2.9}, matrices $M$ and $D_n$ are constant on any edge/face $e\subset\partial K$; for the convection-diffusion-reaction equations \rf{2.11a}-\rf{2.11c}, by a straightforward computation, we have
$$
(M+D_n){\bo u}\cdot\bo w=[2n_1\sqrt{\varepsilon}u_1+2\sqrt{\varepsilon}n_2u_2+(1+\beta\cdot n)u_3]w_3,\;M\bo w\cdot\bo w=w_3^2.
$$
Hence, condition \rf{3.16} holds with $C_M=2\sqrt{\varepsilon}+1+|\beta|_\infty$.
}

For any function $u\in H^1(K)$, the following trace inequality holds.
\begin{equation}
\|u\|^2_{L_2(e)}\leq C\big(h_K^{-1}\|u\|^2_{0,K}+h_K\|\nabla u\|_{0,K}^2\big),\;e\subset\partial K,\,K\in T_h.\la{3.17}
\end{equation}
For a piecewise smooth function $u$, let $u^c$ be its piecewise constant approximation on $T_h$,
$$
u^c(x)=\frac{1}{|K|}\int_Kudx,\;x\in K,\;K\in T_h,
$$
then it holds
\begin{equation}
|u(x)-u^c(x)|\leq Ch_K|u|_{1,\infty,K},\;x\in K,\;K\in T_h.\la{3.16}
\end{equation}
\dl{ 3.2}\quad{\em Assume that $T_h$ is a shape regular partition. Let $\bo u\in [H^{k+1}(\Omega)]^m$ and $\bo u_h\in [V_h]^m$ be the solutions of problems \rf{2.1}-\rf{2.2} and WG equation \rf{2.16}, respectively. Then, we have the following error estimate.}
\begin{equation}
|||Q_h\bo u-\bo u_h|||\leq Ch^{k+\frac{1}{2}}|\bo u|_{k+1},\;k\geq 0.\la{3.18}
\end{equation}
\zm\quad Denote by $\bo e_h=Q_h\bo u-\bo u_h\in [V_h]^m$ the error function. From equation \rf{2.16} and \rf{3.13}, we obtain the error equation for $\bo v\in [V_h]^m$,
\begin{eqnarray}
a(\bo e_h,\bo v)+s(\bo e_h,\bo v)=s(Q_h\bo u,\bo v)-l_1(\bo u,\bo v)-l_2(\bo u,\bo v)-l_3(\bo u,\bo v)-l_4(\bo u,\bo v).\la{3.19}
\end{eqnarray}
Below we estimate the terms $l_i(\bo u,\bo v)$ ($i=1,2,3,4$) and $s(Q_h\bo u,\bo v)$. It follows from the definition of $Q^0$ and the inverse inequality,
$$
(\bo u-Q^0\bo u,\bo A\cdot\nabla \bo v^0)_{T_h}=(\bo u-Q^0\bo u,(\bo A-\bo A^c)\cdot\nabla \bo v^0)_{T_h}\leq C|\bo A|_{1,\infty}h^{k+1}|\bo u|_{k+1}\|\bo v^0\|\,.
$$
Hence, from \rf{3.9} and \rf{3.10}, we obtain
\begin{eqnarray}
|l_1(\bo u,\bo v)+l_2(\bo u,\bo v)|\leq Ch^{k+1}|\bo u|_{k+1}\|\bo v^0\|\leq Ch^{k+1}|\bo u|_{k+1}|||\bo v|||.\la{3.20}
\end{eqnarray}
Next, it follows from the definition of projection operator $Q^b$,
$$
\|\bo u-Q^b\bo u\|_{L_2(e)}^2=\int_e(\bo u-Q^b\bo u)(\bo u-Q^0\bo u)ds\leq \|\bo u-Q^b\bo u\|_{L_2(e)}\|\bo u-Q^0\bo u\|_{L_2(e)},
$$
which implies $\|\bo u-Q^b\bo u\|_{L_2(e)}\leq \|\bo u-Q^0\bo u\|_{L_2(e)},\,e\subset\partial K$. Then, from \rf{3.11} and the trace inequality, we obtain
\begin{eqnarray}
&&|l_3(\bo u,\bo v)|\leq |D_n|_{\infty}\sum_{K\in T_h}\|\bo u-Q^b\bo u\|_{0,\partial K}\|\bo v^0-\bo v^b\|_{0,\partial K}\nonumber\\
&\leq &|D_n|_{\infty}\sum_{K\in T_h}\|\bo u-Q^0\bo u\|_{0,\partial K}\|\bo v^0-\bo v^b\|_{0,\partial K}\leq Ch^{k+\frac{1}{2}}|\bo u|_{k+1}|||\bo v|||.\la{3.21}
\end{eqnarray}
Now, if matrix $M+D_n$ is constant on each element edge $e\subset\partial\Omega$, we have from the definition of $Q^b$ that $l_4(\bo u,\bo v)=0$; otherwise, let $(M+D_n)^c$ be the piecewise constant approximation of matrix $M+D_n$ on element edges, using the trace inequality and inverse inequality, we obtain
\begin{eqnarray}
l_4(\bo u,\bo v)&=&\frac{1}{2}\sum_{e\subset\partial\Omega}\langle(M+D_n-(M+D_n)^c)(\bo u-Q^b\bo u),\bo v^b\rangle_{e}\nonumber\\
&\leq& C\sum_{K\in T_h}h_K|M+D_n|_{1,\partial K,\infty}\|\bo u-Q^b\bo u\|_{0,\partial K}\|\bo v^b\|_{0,\partial K}\nonumber\\
&\leq& C\sum_{K\in T_h}h_K\|\bo u-Q^0\bo u\|_{0,\partial K}\big(\|\bo v^b-\bo v^0\|_{0,\partial K}+\|\bo v^0\|_{0,\partial K}\big)\nonumber\\
&\leq& C\sum_{K\in T_h}h_K\|\bo u-Q^0\bo u\|_{0,\partial K}\big(\|\bo v^b-\bo v^0\|_{0,\partial K}+h_K^{-\frac{1}{2}}\|\bo v^0\|_{0,K}\big)\nonumber\\
&\leq& Ch^{k+1}\|\bo u\|_{k+1}|||\bo v|||.\la{3.22}
\end{eqnarray}
\commentout{
\begin{eqnarray}
&&|l_4(\bo u,\bo v)|\leq \frac{1}{2}C_M\Big(\sum_{e\subset\partial\Omega}\|\bo u-Q^b\bo u\|^2_{L_2(e)}\Big)^{\frac{1}{2}}\langle M\bo v^b,\bo v^b\rangle^{\frac{1}{2}}_{\partial\Omega}\nonumber\\
&\leq& \frac{1}{2}C_M\Big(\sum_{e\subset\partial\Omega}\|\bo u-Q^0\bo u\|^2_{L_2(e)}\Big)^{\frac{1}{2}}\langle M\bo v^b,\bo v^b\rangle^{\frac{1}{2}}_{\partial\Omega}\leq Ch^{k+\frac{1}{2}}|\bo u|_{k+1}|||\bo v|||.\la{3.22}
\end{eqnarray}
}
Finally,
\begin{eqnarray}
&&s(Q_h\bo u,\bo v)=\langle\mu(Q^0\bo u-\bo u+\bo u-Q^b\bo u),\bo v^0-\bo v^b\rangle_{\partial T_h}=\langle\mu(Q^0\bo u-\bo u),\bo v^0-\bo v^b\rangle_{\partial T_h}\nonumber\\
&\leq& \Big(\mu/\mu_0\sum_{K\in T_h}\|Q^0\bo u-\bo u\|_{0,\partial K}^2\Big)^{\frac{1}{2}}
\Big(\mu_0\sum_{K\in T_h}\|\bo v^0-\bo v^b\|_{0,\partial K}^2\Big)^{\frac{1}{2}}\nonumber\\
&\leq& Ch^{k+\frac{1}{2}}|\bo u|_{k+1}|||\bo v|||.\la{3.23}
\end{eqnarray}
Substituting \rf{3.20}--\rf{3.23} into \rf{3.19}, we arrive at
$$
a(\bo e_h,\bo v)+s(\bo e_h,\bo v)\leq Ch^{k+\frac{1}{2}}|\bo u|_{k+1}|||\bo v|||,\;\bo v\in [V_h]^m.
$$
Taking $\bo v=\bo e_h$ and using \rf{3.4a}, it yields
$$
|||\bo e_h|||^2\leq a(\bo e_h,\bo e_h)+s(\bo e_h,\bo e_h)\leq Ch^{k+\frac{1}{2}}|\bo u|_{k+1}|||\bo e_h|||.
$$
The proof is completed.\zb

From Theorem 3.2 and the triangle inequality, we immediately the optimal $L_2$-error estimate,
\begin{equation}
\|\bo u-\bo u_h^0\|\leq Ch^{k+\frac{1}{2}}|\bo u|_{k+1},\;k\geq 0.\la{3.24}
\end{equation}

 Lin et. al. in \cite{Lin} considered a WG method for the singularly perturbed convection-diffusion-reaction problem \rf{2.10}-\rf{2.11}. To compare our WG method with the WG method proposed in \cite{Lin}, we consider the WG method \rf{2.16} for solving the same problem by transforming it into the symmetric hyperbolic systems \rf{2.11a}-\rf{2.11c}. For this problem, since $\rho(D_n)\leq |\beta\cdot n|+\sqrt{\varepsilon}$, we may choose the parameter $\mu=|\beta|_\infty+1$ in WG scheme \rf{2.16}, assuming that $\varepsilon\leq 1$.

Let $\bo u=(u,\sigma_1,\sigma_2)^T$ be the exact solution of problem \rf{2.11a}-\rf{2.11c} and $\bo u_h=(u_h,\sigma_{1,h},\sigma_{2,h})^T\in [V_h]^3$ is the WG solution given by equation \rf{2.16} solving problem \rf{2.11a}-\rf{2.11c}. From Theorem 3.2, we have
\begin{equation}
\|u-u_h^0\|+\|\sigma_1-\sigma_{1,h}^0\|+\|\sigma_{2}-\sigma_{2,h}^0\|\leq Ch^{k+\frac{1}{2}}|\bo u|_{k+1},\;k\geq 0,\la{3.25}
\end{equation}
where constant $C$ is independent of $\varepsilon$. Since $\bo \sigma=-\sqrt{\varepsilon}\,\nabla u$ or $\nabla u=-1/\sqrt{\varepsilon}\,(\sigma_1,\sigma_2)^T$, we may define the weak gradient approximation of $\nabla u$ by setting $\nabla^w u_h=-1/\sqrt{\varepsilon}$ $(\sigma^0_{1,h},\sigma^0_{2,h})^T$. Then, from \rf{3.25}, we obtain the error estimate
\begin{equation}
\|u-u_h^0\|+\sqrt{\varepsilon}\,\|\nabla u-\nabla^wu_h\|\leq Ch^{k+\frac{1}{2}}(|u|_{k+1}+\sqrt{\varepsilon}\,|\nabla u|_{k+1}),\,k\geq 0.\la{3.26}
\end{equation}
On the other hand, Lin et. al. in \cite[Theorem 3.8]{Lin} derived the following error estimate:
\begin{equation}
\|u-u_h^0\|+\sqrt{\varepsilon}\,\|\nabla u-\nabla_wu_h\|\leq Ch^k(h^{\frac{1}{2}}+\sqrt{\varepsilon}\,)\,|u|_{k+1},\,k\geq 0,\la{3.27}
\end{equation}
where $\nabla_wu_h$ is the weak gradient \cite{Wang}. Comparing \rf{3.26} with \rf{3.27}, we see that error estimate \rf{3.26} has half an order higher convergence rate than that given in \rf{3.27}, although the regularity requirement in \rf{3.26} is higher than that in \rf{3.27}. In particular, for the piecewise constant element ($k=0$),
our WG method gives an $O(h^{\frac{1}{2}})$-order convergence rate, but no convergence rate can be obtained from estimate \rf{3.27} if $k=0$.

\section{Numerical experiment}
\setcounter{section}{4} \setcounter{equation}{0}

In this section, we provide some numerical examples to test the performance of the proposed WG method by solving the singularly perturbed  convection-diffusion-reaction equation and the Maxwell's equations.

{\em Singularly perturbed convection-diffusion-reaction problem}

Consider problem:
\begin{eqnarray}
-\varepsilon\triangle u+\boldsymbol{\beta}\cdot\nabla u+\alpha u=f,\;x\in\Omega,\la{4.1}\\
u=0,\;x\in \partial\Omega,\la{4.2}
\end{eqnarray}
where $\Omega\subset R^2$ is a bounded domain, $\varepsilon>0$ is a small parameter and $\alpha-\hbox{div}\beta/2\geq\alpha_0>0$.

We first transform this problem into the symmetric hyperbolic systems \rf{2.11a}-\rf{2.11c} and then use WG method  \rf{2.16} to solve this systems in which the stability parameter $\mu=|\beta|_\infty+1$.
\begin{table}[ht]
  \caption{\label{b-n1} The error profile for solution \eqref{s1} on square grids.}
\begin{center}  \begin{tabular}{c|rr|rr}
\hline level  & $ \|u - u_h^0\|$ &rate &
    $ |||Q_h\bo u - \bo u_h|||$ & rate     \\
\hline
  &\multicolumn{4}{c}{ The $P_1$ WG method, $\varepsilon=10^{-8}$ } \\ \hline
 5&   0.2627E-03 & 1.98&   0.2448E-02 & 1.50\\
 6&   0.6663E-04 & 1.98&   0.8645E-03 & 1.50\\
 7&   0.1681E-04 & 1.99&   0.3053E-03 & 1.50\\
\hline
  &\multicolumn{4}{c}{ The $P_2$ WG method, $\varepsilon=10^{-8}$ } \\ \hline
 5&   0.4300E-03 & 3.05&   0.4185E-02 & 2.55\\
 6&   0.5321E-04 & 3.01&   0.7281E-03 & 2.52\\
 7&   0.6668E-05 & 3.00&   0.1277E-03 & 2.51\\
\hline
  &\multicolumn{4}{c}{ The $P_3$ WG method, $\varepsilon=10^{-8}$ } \\ \hline
 2&   0.9900E-01 & 4.67&   0.3976E+00 & 4.10\\
 3&   0.4772E-02 & 4.37&   0.2765E-01 & 3.85\\
 4&   0.2815E-03 & 4.08&   0.2166E-02 & 3.67\\
\hline
\end{tabular}\end{center} \end{table}

The first example does not have a singularity (boundary layer) and the exact solution of problem \rf{4.1}-\rf{4.2} is
 \begin{align} \la{s1}
  u(x,y)=x(1-x)y(1-y)\quad\hbox{ in } \Omega=(0,1)^2.
 \end{align}
In \eqref{4.1}, set $\beta=(1,2)$ and $\alpha=1$.
We use the uniform square grids in the computation where the first grid consists of one square, and each square is refined into four sub-squares to form
   the next level grid.
The computational results are listed in Table \ref{b-n1} in which $\bo u$ is the solution of the corresponding hyperbolic systems.
Numerical results verify our theoretical analysis.
The convergence is independent of the singular perturbation parameter $\varepsilon$.

\begin{figure}[ht] \begin{center}
\begin{picture}(250,200)(0,0)
 \put(0, 100){\scalebox{0.4}{\includegraphics{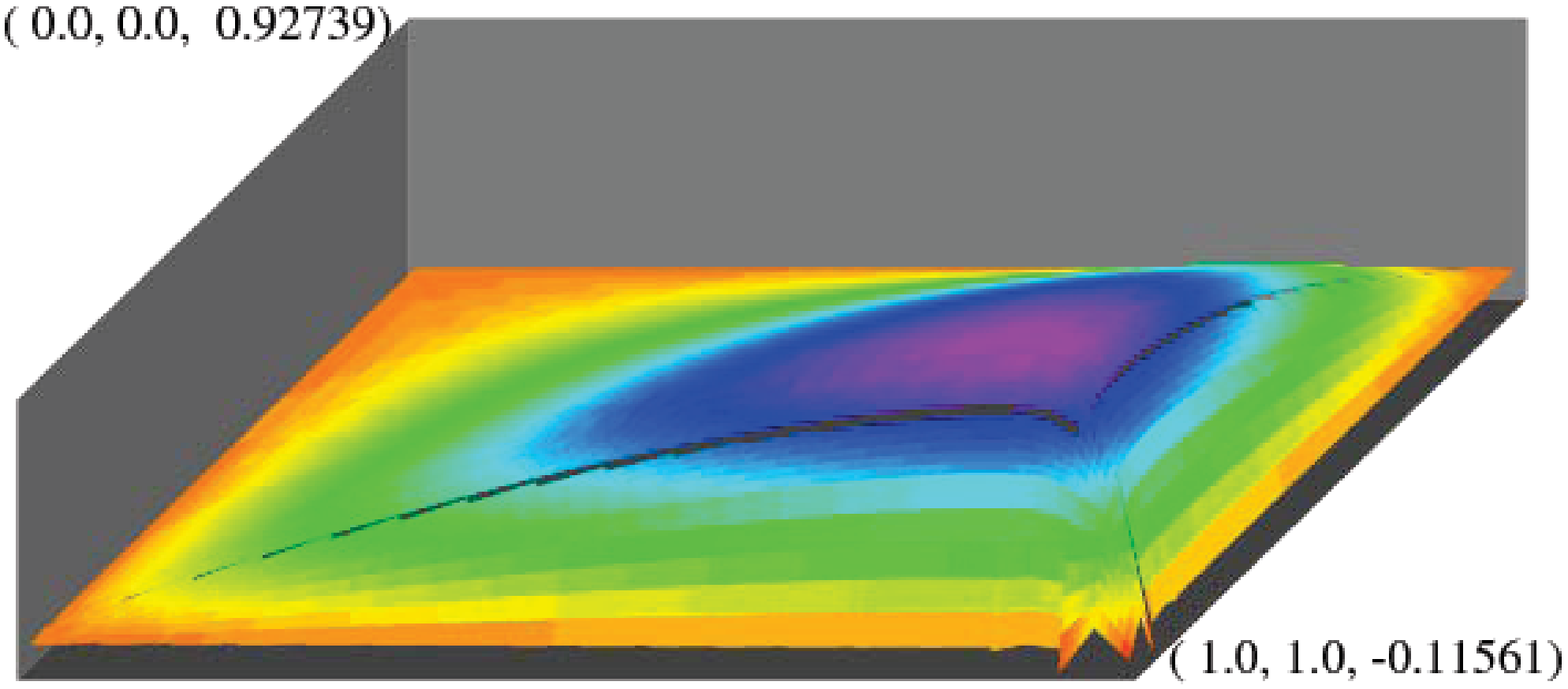}}}
 \put(0,   0){\scalebox{0.4}{\includegraphics{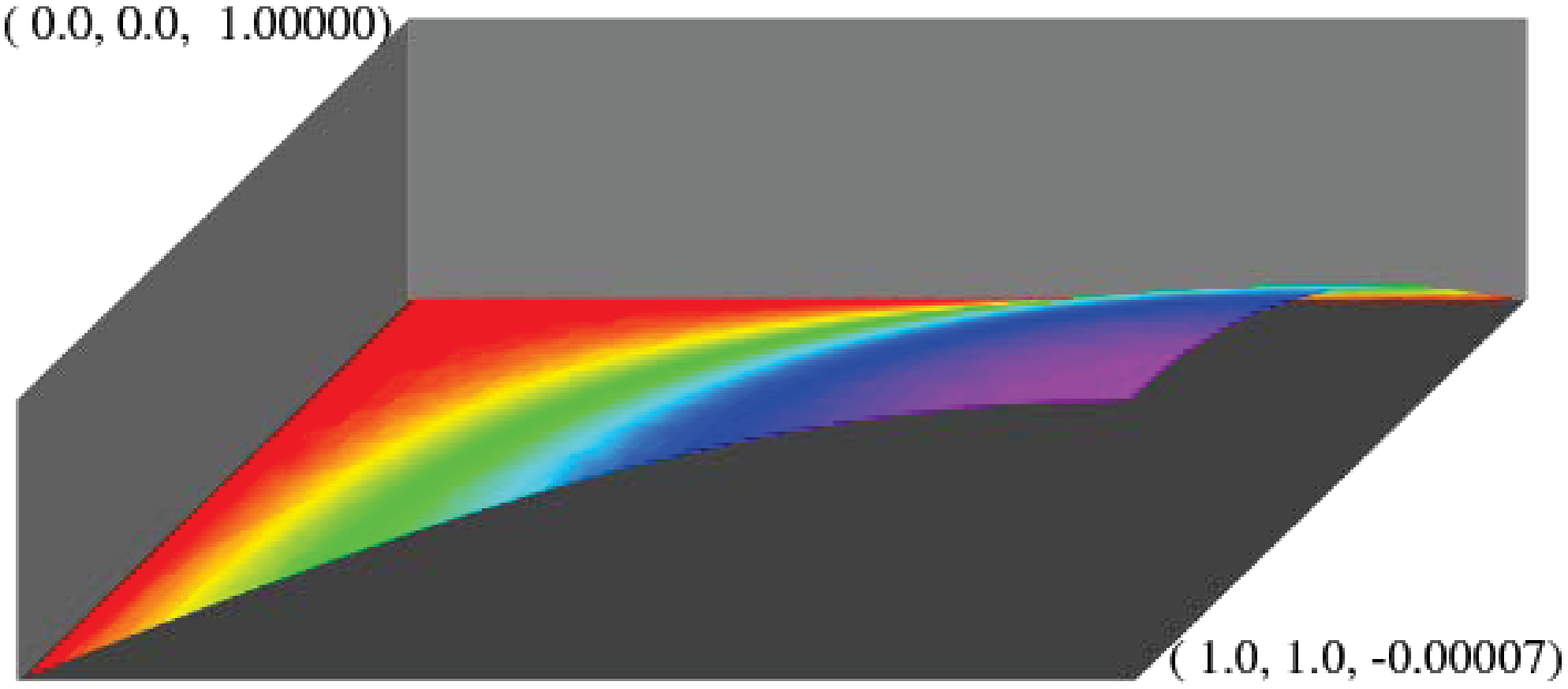}}}
\end{picture}\end{center}
\caption{\label{p1} The $P_2$ WG on the 5th grid
    for solving \rf{4.1}--{4.2}, when $\varepsilon=10^{-3}$ and $10^{-8}$. }
\end{figure}

The second example has a singularity (boundary layer) and the exact solution of problem \eqref{4.1}-\rf{4.2} is
 \begin{align} \la{s2}
  u(x,y)=\sin\frac{\pi x}2 \sin\frac{\pi y}2
     (1-e^{(x-1)/\sqrt{\varepsilon}})(1-e^{(y-1)/\sqrt{\varepsilon}}) \quad\hbox{ in } \Omega=(0,1)^2.
 \end{align}
In \eqref{4.1}, set $\beta=(1,1)$ and $\alpha=1$.
Again the first grid consists of one square, and each square is refined into four sub-squares to form
   the next level grid. The computed solution is plotted in Figure \ref{p1}, where we can see that the boundary layer occurs at the boundary $x=1$ and $y=1$.
The computational results are listed in Table \ref{b-n2}.
They verify our theoretical analysis.
We note that the $P_0$-WG method converges one order higher than the traditional method
   \cite{Lin}, and half an order higher than our theoretic order in the $L_2$-norm.
\begin{table}[ht]
  \caption{\label{b-n2} The error profile for solution \eqref{s2} on square grids.}
\begin{center}  \begin{tabular}{c|rr|rr}
\hline level  & $ \|u - u_h^0\|$ &rate &
    $ |||Q_h\bo u - \bo u_h|||$ & rate     \\ \hline
  &\multicolumn{4}{c}{ The $P_0$ WG method, $\varepsilon=10^{-8}$ } \\
 \hline
 5&   0.1742E-01 & 0.96&   0.3548E-01 & 0.96 \\
 6&   0.8840E-02 & 0.98&   0.1801E-01 & 0.98 \\
 7&   0.4454E-02 & 0.99&   0.9095E-02 & 0.99 \\
\hline
  &\multicolumn{4}{c}{ The $P_1$ WG method, $\varepsilon=10^{-8}$ } \\\hline
 5&   0.4964E-03 & 2.02&   0.6534E-02 & 1.50 \\
 6&   0.1233E-03 & 2.01&   0.2313E-02 & 1.50 \\
 7&   0.3074E-04 & 2.00&   0.8190E-03 & 1.50 \\\hline
  &\multicolumn{4}{c}{ The $P_2$ WG method, $\varepsilon=10^{-8}$ } \\ \hline
 4&   0.1519E-03 & 3.04&   0.1685E-02 & 2.51 \\
 5&   0.1867E-04 & 3.02&   0.2958E-03 & 2.51 \\
 6&   0.2313E-05 & 3.01&   0.5351E-04 & 2.47 \\
\hline
\end{tabular}\end{center} \end{table}

Now, we recompute the solution \eqref{s2} on polygonal grids shown in Figure 2 where
  the central polygon has 12 sides,  and the side ones have 7 sides.
The error and the order of convergence are listed in Table \ref{b-n-n} which
   confirm the $L_2$-error analysis.

\begin{figure}[ht] \label{polygon}\begin{center}
\begin{picture}(330,120)(0,0)
 \put(220, 0){\scalebox{0.18}{\includegraphics{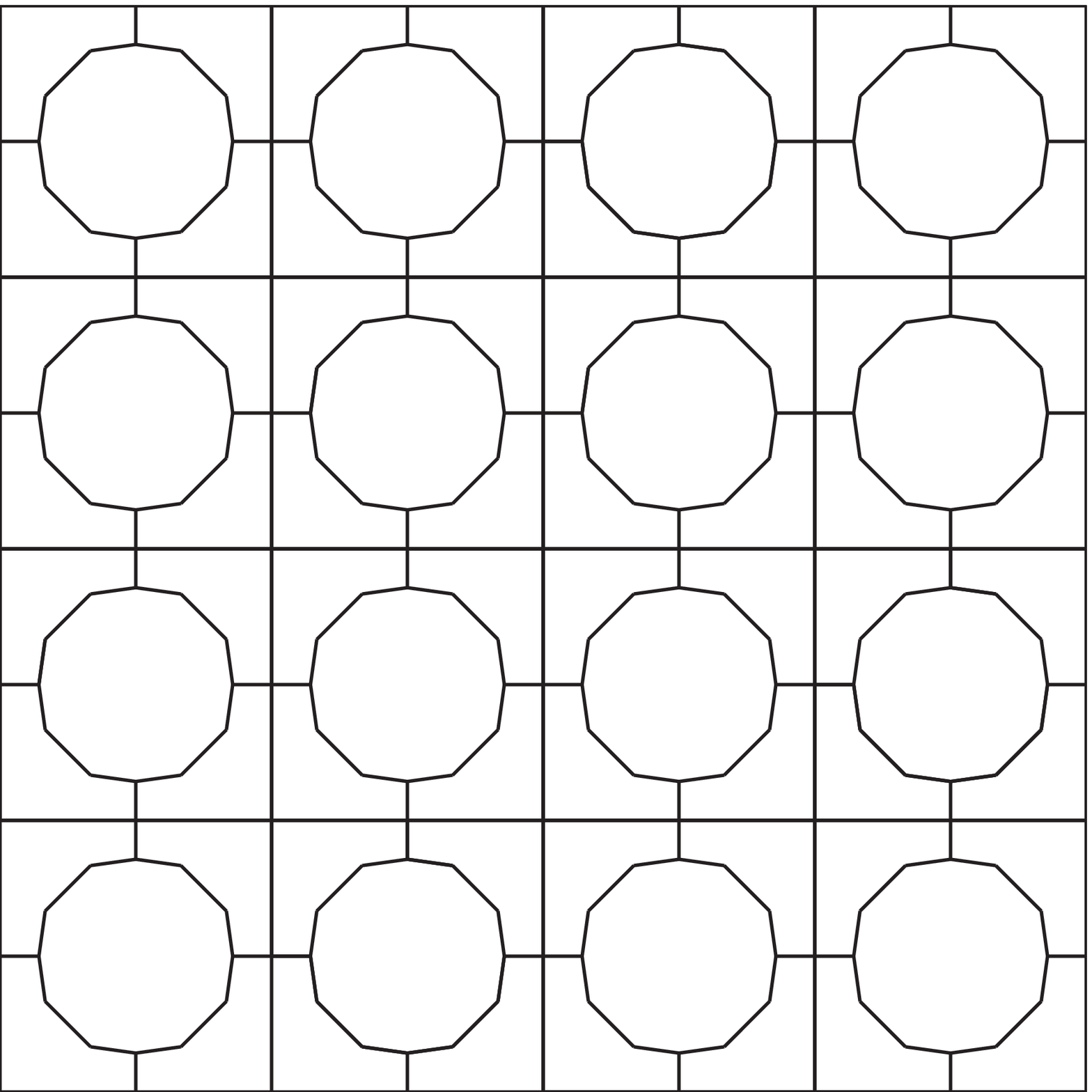}}}
 \put(110, 0){\scalebox{0.18}{\includegraphics{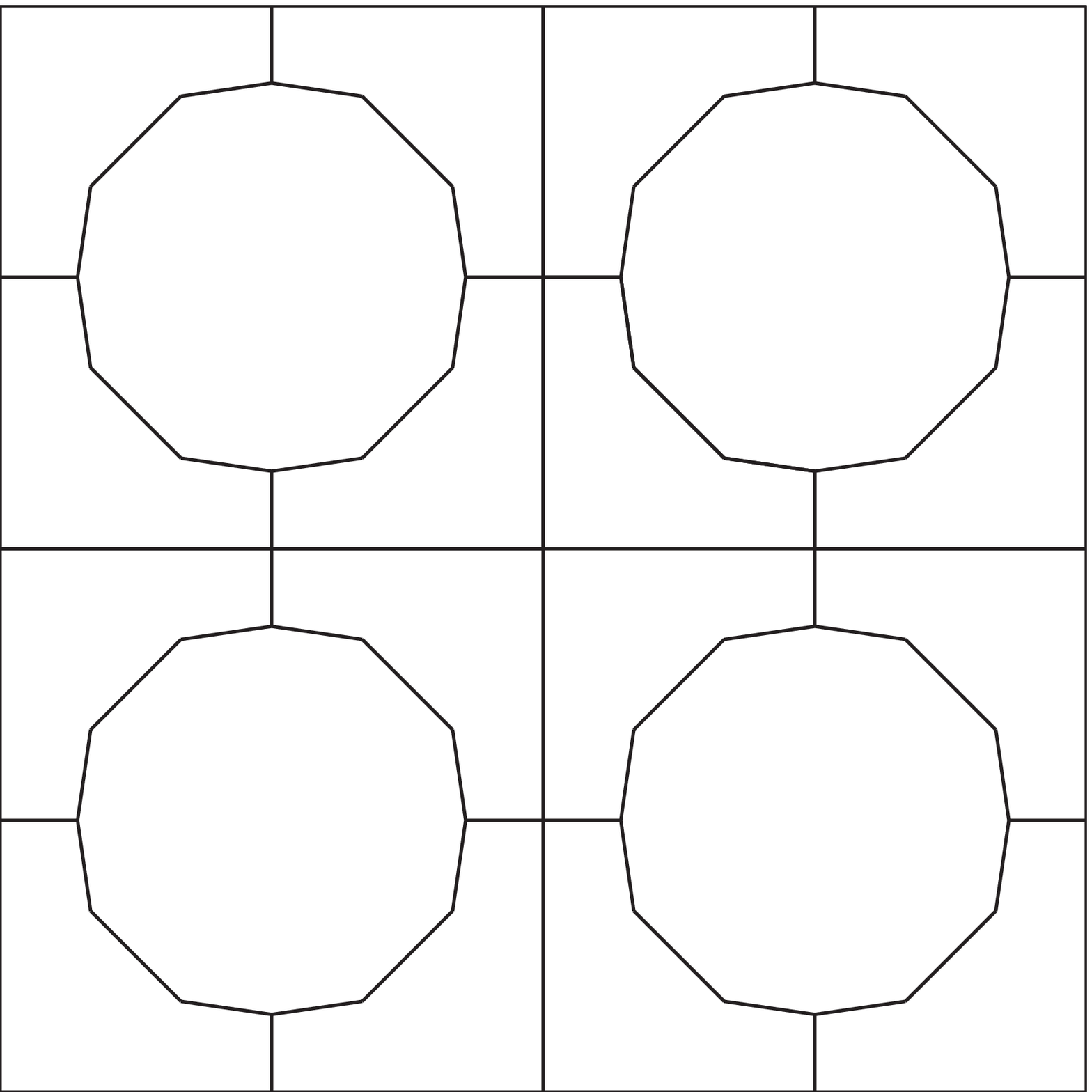}}}
 \put(0,   0){\scalebox{0.18}{\includegraphics{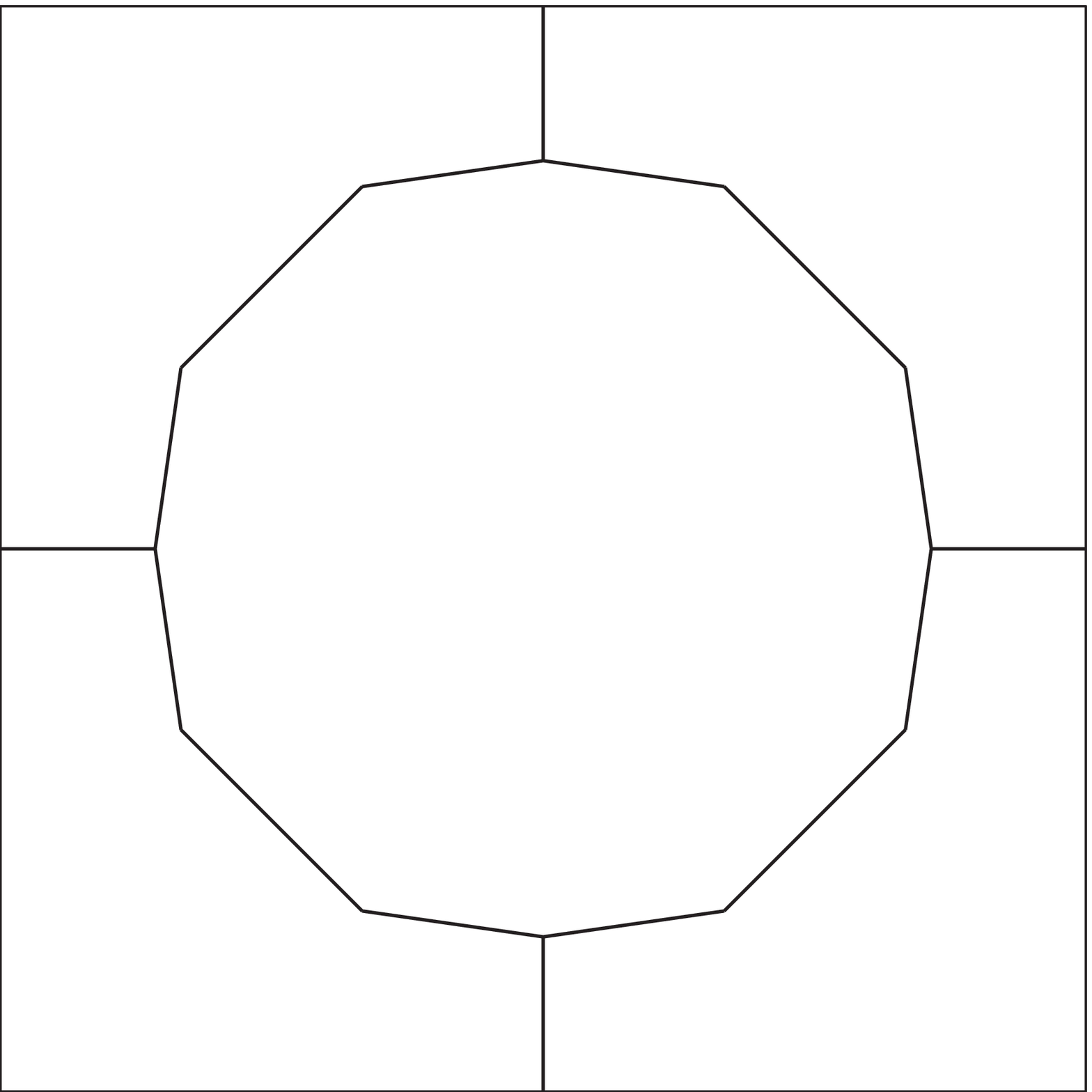}}}
\end{picture}\end{center}
\caption{ The first three polygonal grids
    for solving \eqref{s2}. }
\end{figure}

\begin{table}[ht]
  \caption{\label{b-n-n} The error profile for solving \eqref{s2}
      on polygonal grids (Figure 2).}
\begin{center}  \begin{tabular}{c|rr|rr}
\hline level  & $ \|u - u_h\|_{0}$ &$h^n$ &
    $ |||Q_h u - u_h|||$ & $h^n$     \\ \hline
  &\multicolumn{4}{c}{ The $P_1$ WG method, $\epsilon=10^{-1}$ } \\
 \hline
 4&   0.4714E-03 & 1.97&   0.4274E-02 & 1.47\\
 5&   0.1194E-03 & 1.98&   0.1525E-02 & 1.49\\
 6&   0.3008E-04 & 1.99&   0.5416E-03 & 1.49\\
\hline
  &\multicolumn{4}{c}{ The $P_1$ WG method, $\epsilon=10^{-4}$ } \\
 \hline
 4&   0.1438E+00 & 1.66&   0.1177E+01 & 1.09\\
 5&   0.4464E-01 & 1.69&   0.5384E+00 & 1.13\\
 6&   0.1382E-01 & 1.69&   0.2456E+00 & 1.13\\
\hline
\end{tabular}\end{center} \end{table}

{\em Maxwell's equations}

Consider the Maxwell's equations in a
two-dimensional domain in the following form:
\begin{eqnarray}
\nu H+\nabla\!\wedge E={\bo h},\;(x,y)\in\Omega,\label{4.9}\\
\sigma E-\nabla\!\wedge H=g,\;(x,y)\in\Omega,\label{4.10}\\
E=0,\;(x,y)\in\partial\Omega,\label{4.11}
\end{eqnarray}
where $H=(H_1,H_2)^T$ and $E$ denote the magnetic and electric
fields, respectively, coefficients $\nu>0$ and $\sigma
>0,\,{\bo h}=(h_1,h_2)^T$ and $g$ are known functions, and
$\nabla\!\wedge$ is the $curl$ operator defined by
\begin{eqnarray*}
\nabla\!\wedge E=(\partial_y E\,, -\partial_x
E)^T,\;\;\;\nabla\!\wedge H=\partial_x H_2-\partial_y H_1.
\end{eqnarray*}
Let ${\bo u}=(H_1,H_2,E)^T$. Problem
\rf{4.9}-\rf{4.11} can be written as the positive symmetric hyperbolic system
\begin{eqnarray}
A_1\partial_x{\bo u}+A_2\partial_y{\bo
u}+B{\bo u}={\bo f},\;(x,y)\in\Omega,\label{4.12}\\
 (M-D_n){\bo u}={\bo 0},\;(x,y)\in\partial\Omega,\label{4.13}
\end{eqnarray}
where
$$
A_1=\left(
\begin{array}{ccc}
  0 & 0 & 0 \\
  0 & 0 & -1 \\
  0 & -1 & 0 \\
\end{array}
\right),\;\;
A_2=\left(
\begin{array}{ccc}
  0 & 0 & 1 \\
  0 & 0 & 0 \\
  1 & 0 & 0 \\
\end{array}
\right),\;\; B=\left(
\begin{array}{ccc}
  \nu & 0 & 0 \\
  0 & \nu & 0 \\
  0 & 0 & \sigma \\
\end{array}
\right),\;\; {\bo f}=\left(
\begin{array}{c}
  h_1 \\
  h_2 \\
  g \\
\end{array}
\right),
$$

$$
D_n=\left(
\begin{array}{ccc}
  0 & 0 & n_2 \\
  0 & 0 & -n_1\\
  n_2 & -n_1 & 0 \\
\end{array}
\right),\;\; M=\left(
\begin{array}{ccc}
  0 & 0 & -n_2 \\
  0 & 0 & n_1 \\
  n_2 & -n_1 & 1 \\
\end{array}
\right).
$$
The conditions \rf{2.3}-\rf{2.6} can be verified directly in which $\sigma_0=\min\{\nu,\sigma\}$. Since $\rho(D_n)=1$, we may choose the stability parameter $\mu=1$ in the WG scheme \rf{2.16}.
\begin{table}[ht]
  \caption{\label{b-n3} The error profile for solution \eqref{s3} on square grids.}
\begin{center}  \begin{tabular}{c|rr|rr}
\hline level  & $ \|{\bo u} - {\bo u}_h^0\|$ &rate &
    $ |||{Q_h\bo u} - {\bo u}_h|||$ & rate     \\ \hline
  &\multicolumn{4}{c}{ The $P_1$ WG method  } \\
 \hline
 5&   0.1814E-01 & 1.98&   0.1331E+00 & 1.51\\
 6&   0.4419E-02 & 2.04&   0.4668E-01 & 1.51\\
 7&   0.1067E-02 & 2.05&   0.1642E-01 & 1.51\\
\hline
  &\multicolumn{4}{c}{ The $P_2$ WG method } \\\hline
 5&   0.3386E-02 & 3.04&   0.3361E-01 & 2.55\\
 6&   0.4200E-03 & 3.01&   0.5848E-02 & 2.52\\
 7&   0.5266E-04 & 3.00&   0.1026E-02 & 2.51\\
\hline
  &\multicolumn{4}{c}{ The $P_3$ WG method } \\ \hline
 4&   0.8146E-02 & 4.02&   0.5527E-01 & 3.51\\
 5&   0.5059E-03 & 4.01&   0.4864E-02 & 3.51\\
 6&   0.3154E-04 & 4.00&   0.4290E-03 & 3.50\\
\hline
  &\multicolumn{4}{c}{ The $P_4$ WG method  } \\ \hline
 3&   0.2070E-01 & 5.00&   0.7352E-01 & 4.55\\
 4&   0.6471E-03 & 5.00&   0.3184E-02 & 4.53\\
 5&   0.2022E-04 & 5.00&   0.1393E-03 & 4.52\\
\hline
\end{tabular}\end{center} \end{table}

Now, we apply the WG method \rf{2.16} to solve problem
\rf{4.9}-\rf{4.11}. Let domain $\Omega=(0,1)\times(0,1)$ and the
coefficients are $\nu=1$ and $\sigma=1$. We take the exact
solution
\begin{align}\la{s3}
{\bo u}=\begin{pmatrix} H_1\\ H_2 \\ E\end{pmatrix}
    =\begin{pmatrix} -\partial E/\partial y\\ \partial E/\partial x\\8x(1-x)y(1-y)\end{pmatrix}.
\end{align}
We compute the problem on uniform square grids where the first grid consists of one square, the domain, and each square is refined into four sub-squares to form the next level grid.
The convergence results are listed in Table \ref{b-n3}, which match perfectly our theoretic order in the $|||\cdot|||$-norm and are half an order higher than our theoretic order in the $L_2$-norm.

\section{Conclusion}
\setcounter{section}{5} \setcounter{equation}{0}
We present and analyze a weak Galerkin finite element (WG) method for solving the symmetric hyperbolic systems. This method is highly flexible by allowing to use the discontinuous finite elements on element and its boundary independently of each other. We establish a stable weak Galerkin scheme and derive the optimal $L_2$-error estimate of order $O(h^{k+\frac{1}{2}})$ for $k\geq 0$. This WG method can be applied to many important physical problems. For the singularly perturbed convection-diffusion-reaction equation, we derive an $\varepsilon$-uniform error estimate of order $k+1/2$. Numerical examples show the effectiveness of the proposed WG method.

\section*{Acknowledgments}
This work was supported by the State Key Laboratory of Synthetical
Automation for Process Industries Fundamental Research Funds, No. 2013ZCX02.

\end{document}